# Filtering theory for a weakly coloured noise process

*Shaival H. Nagarsheth (Corresponding author)*
*Electrical Engineering Department*
*Sardar Vallabhbhai National Institute of Technology, Surat, India.*
*Email: shn411@gmail.com, Mob. No.: +91-9727161696, ORCID ID: 0000-0001-9867-8167*

*Dhruvi S. Bhatt*
*Electrical Engineering Department*
*Sardar Vallabhbhai National Institute of Technology, Surat, India.*
*Email: bhattdhruvi427@gmail.com, ORCID ID: 0000-0002-4083-6422*

*Shambhu N. Sharma*
*Electrical Engineering Department*
*Sardar Vallabhbhai National Institute of Technology, Surat, India*
*Email: snsvolterra@gmail.com*

**Abstract**

The problem of analyzing the Itô stochastic differential system and its filtering has received attention. The classical approach to accomplish filtering for the Itô SDE is the Kushner equation. In contrast to the classical filtering approach, this paper presents filtering for the stochastic differential system affected by weakly coloured noise, i.e. $\dot{x}_t = f(x_t) + g(x_t)\xi_t$, where the input process $\xi_t$ is a weakly coloured noise process. As a special case, the process $\xi_t$ can be regarded as the Ornstein-Uhlenbeck (OU) process, i.e. $d\xi_t = -\alpha \xi_t dt + \beta dB_t$, where $\alpha > 0$. More precisely, the filtering model of this paper can be cast as

$$\dot{x}_t = f(x_t) + g(x_t)\xi_t,$$

$$z_t = \int_{t_0}^{t} h(x_\tau) d\tau + \eta_t,$$

where $h(x_t)$ is the measurement non-linearity and $\eta = \{\eta_t, t_0 \leq t < \infty\}$ is the Brownian motion process. The former expression describes the structure of a noisy dynamical system and the latter is observation equation. The theory of this paper is based on a pioneering contribution of Stratonovich involving the perturbation-theoretic approach to noisy dynamical systems in combination with the notion of the 'filtering density' evolution. Making the use of the filtering density evolution equation, the stochastic evolution of condition moment is derived. A scalar Duffing system driven by the OU process is employed to test the effectiveness of the filtering theory of the paper. Numerical simulations involving four different sets of initial conditions and system parameters are utilized to examine the efficacy of the filtering algorithm of this paper.

**Key words**: non-Markovian stochastic system, Ornstein-Uhlenbeck process, Brownian motion process, classical filtering equations, Kushner equation, filtering density, stochastic differential equation.

## 1. Introduction

The Stochastic Differential Equation (SDE) formalism has found remarkable success in diverse field, e.g. adaptive control, satellite trajectory estimations, helicopter rotor, stochastic networks, mathematical finance, neuronal activity, protein kinematics (Kloeden and Platen 1991, Campen 2007). The Itô stochastic differential equation is widely used in analyzing randomly perturbed dynamical systems. Rigorous mathematical treatments about stochastic differential equations can be found in Karatzas and Shreve (1991), Strook and Varadhan (1979). In 1945, Kiyoshi Itô developed the rigorous mathematical framework for the Brownian motion process as well as the non-linear stochastic differential system driven by the Brownian motion. The Itô contribution is popularly known as Itô calculus. The white noise can be regarded as an informal non-existent time derivative $\dot{B}_t$ of the Brownian motion process. Furthermore, in 1975, Hida developed white noise calculus in which white noise can be regarded as the multiplication operator (Hida *et al.* 1993; Kuo 2009). Filtering for the stochastic differential system in Itô sense using the Kushner's equation has been studied extensively in the literature. The Kushner's equation, a filtering density evolution equation, assumes the structure of a partial-integro stochastic differential equation. Liptser and Shiryayev (1977) developed the stochastic evolution of filtering density using a different method and different notations. In their proof of the filtering density, the stochastic evolution of the conditional moment was derived. The stochastic evolution of the conditional characteristic function becomes a special case of the conditional moment evolution. Interestingly, the filtering density can be defined as the inverse Fourier transform of the conditional characteristic function, which leads to the filtering density evolution equation. Stratonovich developed the filtering density evolution equation for the stochastic differential equation with $\frac{1}{2}$-differential. The stochastic differential equation of the form $\dot{x}_t = f(x_t) + g(x_t)\xi_t$ is a structure, where the input process $\xi_t$ is arbitrary. Here, the arbitrary process means no restrictions imposed on the stochastic process as well as generates a non-Markov process. For $\xi_t = \dot{B}_t$ in combination with the term $dB_t = \dot{B}_t dt$, the SDE reduces to the Itô formalism. The problem of analyzing the arbitrary process-driven dynamical system can be accomplished using the Stratonovich perturbation-theoretic approach as well as the functional calculus approach. The functional calculus approach allows deriving the master equation for coloured noise. The master equation for coloured noise is an infinite series. After introducing the notion of zero mean, stationarity, and the Gaussian assumption for the input process $\xi_t$, the master equation assumes a quite simplified structure involving the functional derivative $\frac{\delta x_t}{\delta \xi_s}$ within the single integral sign, where the process $x_t$ is non-Markov, see Hänggi (1995, p. 85). In physical problems, the random perturbation in dynamical systems has non-zero correlation time, e.g. the Ornstein-Uhlenbeck process. The Ornstein-Uhlenbeck process has an interpretation as a coloured noise process in dynamical systems theory. Bhatt and Karandikar (1999) have developed a filtering equation for the filtering model involving the OU process as the observation noise, see Mandrekar and Mandal (2000) as well. Their filtering model involves the structure of observations, i.e. $z_t = \int_{t_0}^{t} \int_{t_0 - \frac{1}{\alpha}}^{\tau} h(x_s, s) ds + N_t^{\beta}$, where $N_t^{\beta}$ is the displacement OU process, $\alpha > 0$, and the damping ratio $\beta > 0$. The Kushner theory and the Itô SDE, can be found in Jazwinski (1970, p. 178), Liptser and Shiryayev (1977, p. 318). Looking into the reality of the physical situation for the dynamic systems, affected by weakly coloured noise, a filtering theory in the Kushner setting merits investigations for dynamic systems driven by the weakly coloured noise.

This paper is devoted to develop and analyze filtering theory for weakly coloured noise driven stochastic differential system, a special case of the non-Markovian SDE. Note that the 'weakly coloured noise' process considered in this paper is a zero mean, stationary process with finite, non-zero, smaller correlation time. The filtering model of this paper can be regarded as an *extension* of the filtering model for the Itô stochastic differential system. The classical filtering equations become a special case of the filtering equations of this paper. The theory of this paper is grounded on the fact that an actual output process can be replaced with a Markov process, where the input process has a smaller correlation time. Subsequently, a Kushner-type equation is derived, which allows obtaining the stochastic evolution of conditional moment. The developed filtering theory of the paper is applied to a scalar Duffing system driven by the OU process. Numerical simulations are carried out under a variety of conditions to test the efficacy and usefulness of the developed filtering theory.

This paper is organized as follows: in section (2), we derive the stochastic evolutions of conditional mean and variance. Section (3) discusses the usefulness of the theory of the paper by analyzing the scalar damped Duffing SDE. Section (3) is about numerical simulations based on the filtering equations of this paper. Section (4) encompasses concluding remarks. *Appendix* A of the paper discusses an alternative proof of the Fokker-Planck equation for the nonlinear stochastic differential system driven by an OU process.

## 2. Filtering equations for non-Markovian stochastic differential systems

Filtering for the Itô stochastic differential system is accomplished using the filtering model (Jazwinski 1970), i.e.
$$dx_t = f(x_t, t)dt + G(x_t, t)dB_t,$$
and
$$dz_t = h(x_t, t)dt + d\eta_t,$$
where $G(x_t, t)$ is the process noise coefficient matrix and var $(dB_t) = Idt$. In contrast to the Itô stochastic differential system, this paper considers a non-Markovian stochastic differential system, i.e. $\dot{x}_t = F(x_t, \xi_t)$. As a special case, the stochastic differential system
$$\dot{x}_t = f(x_t) + g(x_t)\xi_t \tag{1}$$
with the observation equation

$$z_t = \int_{t_0}^{t} h(x_\tau) d\tau + \eta_t, \qquad (2)$$

be the subject of investigation, where var $(d\eta_t) = \varphi_\eta dt$,. The evolution of conditional probability density for given initial states for a non-Markov process assumes the structure of an infinite series, popularly known as stochastic equation (kinetic equation). A proof of the kinetic equation begins from the relationship $p(x_1) = \int p(x_1|x_2)p(x_2)dx_2$, where $x_{t_1} = x_1$, $x_{t_2} = x$, and $t_1 = t+\tau, t_2 = t$. The relationship between the transition probability density and the conditional characteristic function using the inverse Fourier integral and some simplifications lead to the stochastic equation. The stochastic equation can be stated as $\dot{p}(x) = \sum_{1 \leq n} \frac{1}{n!}(-\frac{\partial}{\partial x})^n k_n(x)p(x)$, where $k_n(x) = \frac{\langle x_{t+\tau} - x \rangle}{\tau}, \tau \to 0$. Most notably, the stochastic equation for the weakly coloured noise driven stochastic differential equation involves only two spatial derivatives and the resulting equation is the Fokker-Planck equation, i.e.

$\dot{p}(x) = -\frac{\partial}{\partial x} k_1(x)p(x) + \frac{1}{2}\frac{\partial^2 k_2(x)p(x)}{\partial x^2}$. The higher-order terms involving $k_n(x), n > 2$ will vanish. The Fokker-Planck equation describes the evolution of conditional probability density for a Markov process with given initial states, which satisfies a stationary Gaussian white noise driven stochastic differential equation. Thus, the weakly coloured noise driven and the white noise driven stochastic differential systems would be associated with the same Fokker-Planck equation and both stochastic differential systems can be considered 'stochastically' equivalent. A simple calculation shows that the white noise driven stochastic differential equation associated with the Fokker-Planck equation (Stratonovich 1963, p. 97) assumes the structure

$$\dot{x}_t = (k_1(x) - \frac{k_2'(x)}{4}) + \sqrt{k_2(x)} w_0(t), \qquad (3)$$

where the input process $w_0(t)$ is a zero mean, stationary, Gaussian white noise process. The coefficients $k_1(x)$ and $k_2(x)$ for the stochastic differential system $\dot{x}_t = f(x_t) + g(x_t)\xi_t$, can be re-cast as (Stratonovich 1963)

$$k_1(x) = f + \frac{\mu_1}{2} gg' + \mu_2 g'(\frac{f}{g})', \qquad (4)$$

$$k_2(x) = \mu_1 g^2 + 2\mu_2 g^3 (\frac{f}{g})', \qquad (5)$$

where $\mu_1 = 2\int_{-\infty}^{0} R_{\xi\xi}(\tau)d\tau$, $\mu_2 = \int_{-\infty}^{0} |\tau| R_{\xi\xi}(\tau)d\tau$, and $R_{\xi\xi}(\tau) = E\xi_{t+\tau}\xi_t$. The terms $g'$ and $(\frac{f}{g})'$ denote the derivatives with respect to the spatial variable $x$. Note that equations (4)-(5) of this paper are valid for the weakly coloured noise input process $\xi_t$. Equation (3) in conjunction with equations (4)-(5) assumes the structure

$$\dot{x}_t = f(x_t) - \frac{1}{2}\mu_2 g^2(x_t) g'(x_t) (\frac{f(x_t)}{g(x_t)})' - \frac{\mu_2}{2} g^3(x_t)(\frac{f(x_t)}{g(x_t)})''$$
$$+ \sqrt{\mu_1} g(x_t) \sqrt{1 + \frac{2\mu_2}{\mu_1} g(x_t)(\frac{f(x_t)}{g(x_t)})'} w_0(t). \qquad (6)$$

Note that equation (6) is *stochastically* equivalent to the SDE, i.e. $\dot{x}_t = f(x_t) + g(x_t)\xi_t$, since both are associated with the same Fokker-Planck equation. Equation (6) can be reformulated in the Itô sense, i.e.

$$dx_t = (f(x_t) - \frac{1}{2}\mu_2 g^2(x_t) g'(x_t)(\frac{f(x_t)}{g(x_t)})' - \frac{\mu_2}{2} g^3(x_t)(\frac{f(x_t)}{g(x_t)})'')dt$$
$$+ \sqrt{\mu_1} g(x_t)\sqrt{1 + \frac{2\mu_2}{\mu_1} g(x_t)(\frac{f(x_t)}{g(x_t)})'} dB_t, \qquad (7)$$

where

$$a(x_t) = f(x_t) - \frac{1}{2}\mu_2 g^2(x_t) g'(x_t)(\frac{f(x_t)}{g(x_t)})' - \frac{\mu_2}{2} g^3(x_t)(\frac{f(x_t)}{g(x_t)})'', \qquad (8)$$

$$b(x_t) = \sqrt{\mu_1} g(x_t) \sqrt{1 + \frac{2\mu_2}{\mu_1} g(x_t)(\frac{f(x_t)}{g(x_t)})'}. \qquad (9)$$

Equation (7), the Ito stochastic differential equation, describes a rigorous interpretation of the SDE in contrast to equation (6).

### 2.1 The OU process-driven stochastic differential system

The OU process is a Gauss-Markov process and satisfies the stochastic differential equation $d\xi_t = -\frac{1}{\tau_{cor}} \xi_t dt + \frac{\sqrt{2D}}{\tau_{cor}} dB(t)$, where the term $\frac{\sqrt{2D}}{\tau_{cor}}$ has an interpretation as the process noise coefficient. For the OU process, a weakly coloured noise process, the terms $\mu_1, \mu_2$ of equations (4)-(5) become

$$\mu_1 = 2D, \quad \mu_2 = D\tau_{cor}. \qquad (10)$$

From equations (8)-(9), we have the following system non-linearity $a(x_t)$ and process noise coefficient $b(x_t)$ for the OU process-driven stochastic differential system:

$$a(x_t) = f(x_t) - \frac{D\tau_{cor}}{2} g^2(x_t) g'(x_t)(\frac{f(x_t)}{g(x_t)})' - \frac{D\tau_{cor}}{2} g^3(x_t)(\frac{f(x_t)}{g(x_t)})'', \qquad (11)$$

$$b(x_t) = \sqrt{2D} g(x_t) \sqrt{1 + \tau_{cor} g(x_t)(\frac{f(x_t)}{g(x_t)})'}. \qquad (12)$$

This paper is intended for developing the filtering equations of the weakly coloured noise driven SDE in which the structure of the observation equation is driven by the Brownian noise. This structure of the observation equation allows us to write the filtering density $p(x, t|z_\tau, t_0 \leq \tau \leq t)$ in the Kushner setting (Jazwinski 1970, p. 178). The Kushner-type equation for filtering of the system of equations (1)-(2), where equation (1) is 'stochastically equivalent' to equation (7), can be stated as

$$dp = (-\frac{\partial f(x)p}{\partial x} + \frac{\mu_2}{2} \frac{\partial g^2(x)g'(x)(\frac{f(x)}{g(x)})'p}{\partial x} + \frac{\mu_2}{2} \frac{\partial g^3(x)(\frac{f(x)}{g(x)})''p}{\partial x}$$

$$+ \frac{\mu_1}{2} \frac{\partial^2 g^2(x)(1 + 2\frac{\mu_2}{\mu_1} g(x)(\frac{f(x)}{g(x)})')p}{\partial x^2})dt + (h(x) - \hat{h}(x_t))\varphi_\eta^{-1}(dz_t - \hat{h}(x_t)dt)p,$$

$$= L(p)dt + (h(x) - \hat{h}(x_t))\varphi_\eta^{-1}(dz_t - \hat{h}(x_t)dt)p, \qquad (13)$$

where $p = p(x, t|z_\tau, t_0 \leq \tau \leq t)$. From equation (13), the Fokker-Planck operator for a stochastic differential system of equation (2) can be stated as

$$L(.) = -\frac{\partial f(x)(.)}{\partial x} + \frac{\mu_2}{2} \frac{\partial g^2(x)g'(x)(\frac{f(x)}{g(x)})'(.)}{\partial x} + \frac{\mu_2}{2} \frac{\partial g^3(x)(\frac{f(x)}{g(x)})''(.)}{\partial x}$$

$$+ \frac{\mu_1}{2} \frac{\partial^2 g^2(x)(1 + 2\frac{\mu_2}{\mu_1} g(x)(\frac{f(x)}{g(x)})')(.)}{\partial x^2}. \qquad (14)$$

Equations (13)-(14) in combination with equation (10) lead to the Kushner-type equation for the OU process-driven stochastic differential system, i.e.

$$dp = (-\frac{\partial f(x)p}{\partial x} + \frac{D\tau_{cor}}{2} \frac{\partial g^2(x)g'(x)(\frac{f(x)}{g(x)})p}{\partial x} + \frac{D\tau_{cor}}{2} \frac{\partial g^3(x)(\frac{f(x)}{g(x)})''p}{\partial x}$$

$$+ D \frac{\partial^2 g^2(x)(1 + 2\frac{\mu_2}{\mu_1} g(x)(\frac{f(x)}{g(x)})')p}{\partial x^2})dt$$

$$+ (h(x) - \hat{h}(x_t))\varphi_\eta^{-1}(dz_t - \hat{h}(x_t)dt)p.$$

## 2.2 Stochastic evolution of conditional moment

Consider the scalar function $\phi(x_t)$, system nonlinearity $f(x_t)$, and the process noise coefficient $g(x_t)$ are twice continuously differentiable with respect to $x_t$. More precisely, $\phi: R \to R$, $R$ is the solution space of the equation (1). The stochastic evolution $d\hat{\phi}(x_t)$ of conditional moment is defined as

$$d\hat{\phi}(x_t) = \int \phi(x) dp(x, t|z_\tau, t_0 \leq \tau \leq t) dx, \tag{15}$$

where $\hat{\phi}(x_t) = E(\phi(x_t)|z_\tau, t_0 \leq \tau \leq t)$. From equations (13)-(14), we have

$$d\hat{\phi}(x_t) = \int \phi(x)(L(p)dt + (h(x) - \hat{h}(x_t))\varphi_\eta^{-1}(dz_t - \hat{h}(x_t)dt)p)dx.$$

$$= \langle \phi, Lp \rangle dt + (\widehat{\phi h} - \hat{\varphi}\hat{h})\varphi_n^{-1}(dz_t - \hat{h}dt)$$

$$= \langle L'\phi, p \rangle dt + (\widehat{\phi h} - \hat{\varphi}\hat{h})\varphi_n^{-1}(dz_t - \hat{h}dt), \tag{16}$$

where the backward operator (Karatzas and Shreve 1991, p. 281) for the stochastic differential system of equation (2) is

$$L'(.) = f \frac{\partial(.)}{\partial x} - \frac{\mu_2}{2} g^2(x)g'(x)(\frac{f(x)}{g(x)})' \frac{\partial(.)}{\partial x} - \frac{\mu_2}{2} g^3(x)(\frac{f(x)}{g(x)})'' \frac{\partial(.)}{\partial x}$$

$$+ \frac{\mu_1}{2} g^2(x)(1 + 2\frac{\mu_2}{\mu_1} g(x)(\frac{f(x)}{g(x)})') \frac{\partial^2(.)}{\partial x^2}, \tag{17}$$

From equations (16)-(17), we have

$$d\hat{\phi}(x_t) = (\overbrace{(f(x_t) - \frac{1}{2}\mu_2 g^2(x_t)g'(x_t)(\frac{f(x_t)}{g(x_t)})' - \frac{\mu_2}{2} g^3(x_t)(\frac{f(x_t)}{g(x_t)})'')\frac{\partial \phi(x_t)}{\partial x_t}}$$

$$+ \frac{1}{2}\overbrace{\mu_1 g^2(x_t)(1 + 2\frac{\mu_2}{\mu_1} g(x_t)(\frac{f(x_t)}{g(x_t)})')\frac{\partial^2 \phi(x_t)}{\partial x_t^2}})dt + (\widehat{\phi h} - \hat{\phi}\hat{h})\varphi_\eta^{-1}(dz_t - \hat{h}dt). \tag{18}$$

where $\widehat{\phi h} = E(\phi(x_t)h(x_t)|z_\tau, t_0 \leq \tau \leq t)$. Interestingly, equation (18) is a consequence of the notion of conditional expectation, integration by part formula, and adjoint property of the Fokker-Planck operator. For the correlation time $\mu_2 = \tau_{cor} = 0$, we arrive at standard filtering equations, i.e. the scalar version of Lemma of (6.3) of Jazwinski (1970, p. 183). For $\phi(x_t) = x_t$ and $\phi(x_t) = \tilde{x}_t^2$, we have the stochastic evolutions of conditional mean and variance, where $\tilde{x}_t = x_t - \hat{x}_t$,

$$d\hat{x}_t = \overbrace{(f(x_t) - \frac{1}{2}c_2 g^2(x_t)g'(x_t)(\frac{f(x_t)}{g(x_t)})' - \frac{c_2}{2} g^3(x_t)(\frac{f(x_t)}{g(x_t)})'')dt},$$

$$+ (\widehat{\phi h} - \hat{\phi}\hat{h})\varphi_\eta^{-1}(dz_t - \hat{h}dt). \tag{19}$$

$$dP_t = \overbrace{(2x_t(f(x_t) - \frac{1}{2}c_2 g^2(x_t)g'(x_t)(\frac{f(x_t)}{g(x_t)})' - \frac{c_2}{2} g^3(x_t)(\frac{f(x_t)}{g(x_t)})'')}$$

$$- 2\hat{x}_t \overbrace{(f(x_t) - \frac{1}{2}c_2 g^2(x_t)g'(x_t)(\frac{f(x_t)}{g(x_t)})' - \frac{c_2}{2} g^3(x_t)(\frac{f(x_t)}{g(x_t)})'')}$$

$$+ \overbrace{(c_1 g^2(x_t)(1 + 2\frac{c_2}{c_1} g(x_t)(\frac{f(x_t)}{g(x_t)})'))dt} + (\widehat{\phi h} - \hat{\phi}\hat{h})\varphi_\eta^{-1}(dz_t - \hat{h}dt). \tag{20}$$

The coefficients associated with the term '$dt$' of equations (18)-(20) can be regarded as the conditional expectation of the terms within the brackets. It is interesting to note that the above exact filtering equations (19)-(20), can be regarded as an extension of the scalar version of Theorem (6.6) stated in Jazwinski (1970, p.184).

The system non-linearity $f(x_t)$, the process noise coefficient term $g(x_t)$, and their spatial derivatives are associated with the weakly coloured noise driven stochastic differential system of equation (1). Note that the stochastic *filtering equations* comprise the stochastic evolutions of conditional mean and variance. The analytical and numerical solutions to the exact filtering equations are intractable since they become infinite dimensional and require the knowledge of higher-order moments (Kushner 1967). For these reasons, this paper applies the second-order approximation to the system non-linearity, measurement non-linearity, the diffusion coefficient as well as nearly Gaussian assumptions. The second-order filtering equations offer a simplified analysis as well as preserve some of the qualitative characteristics of the exact filtering equations. The second-order filtering equations for the stochastic differential system described by equation (1) can be derived by introducing the second-order partials of the terms $a(x_t)$, $b(x_t)$ of equations (8)-(9) and the measurement non-linearity $h(x_t)$, which are evaluated at the conditional mean of the state $x_t$, into equations (19)-(20). After a few steps of simplifications, we arrive at

$$d\hat{x}_t = (f(\hat{x}_t) - \frac{\mu_2}{2} g^2(\hat{x}_t) g'(\hat{x}_t)(\frac{f(\hat{x}_t)}{g(\hat{x}_t)})' - \frac{\mu_2}{2} g^3(\hat{x}_t)(\frac{f(\hat{x}_t)}{g(\hat{x}_t)})''$$
$$+ \frac{1}{2} P_t (f''(\hat{x}_t) - \frac{\mu_2}{2}(g^3(\hat{x}_t)\frac{d^4}{d\hat{x}_t^4}(\frac{f(\hat{x}_t)}{g(\hat{x}_t)}) + 7g^2(\hat{x}_t)g'(\hat{x}_t)(\frac{f(\hat{x}_t)}{g(\hat{x}_t)})'''$$
$$+ 5g^2(\hat{x}_t)g''(\hat{x}_t)(\frac{f(\hat{x}_t)}{g(\hat{x}_t)})'' + 10g(\hat{x}_t)(g'(\hat{x}_t))^2(\frac{f(\hat{x}_t)}{g(\hat{x}_t)})''$$
$$+ g^2(\hat{x}_t)g'''(\hat{x}_t)(\frac{f(\hat{x}_t)}{g(\hat{x}_t)})' + 2(g'(\hat{x}_t))^3(\frac{f(\hat{x}_t)}{g(\hat{x}_t)})'$$
$$+ 6g(\hat{x}_t)g'(\hat{x}_t)g''(\hat{x}_t)(\frac{f(\hat{x}_t)}{g(\hat{x}_t)})')))dt$$
$$+ P_t h'(\hat{x}_t)\varphi_\eta^{-1}(dz_t - (h(\hat{x}_t) + \frac{1}{2}P_t h''(\hat{x}_t))dt), \qquad (21)$$

$$dP_t = (2P_t(f'(\hat{x}_t) - \frac{1}{2}\mu_2(4g^2(\hat{x}_t)g'(\hat{x}_t)(\frac{f(\hat{x}_t)}{g(\hat{x}_t)})'' + g^2(\hat{x}_t)g''(\hat{x}_t)(\frac{f(\hat{x}_t)}{g(\hat{x}_t)})'$$
$$+ 2g(\hat{x}_t)(g'(\hat{x}_t))^2(\frac{f(\hat{x}_t)}{g(\hat{x}_t)})' + g^3(\hat{x}_t)(\frac{f(\hat{x}_t)}{g(\hat{x}_t)})'''))$$
$$+ \mu_1 g^2(\hat{x}_t)(1 + 2\frac{\mu_2}{\mu_1}g(\hat{x}_t)(\frac{f(\hat{x}_t)}{g(\hat{x}_t)})')$$
$$+ P_t(\mu_2 g^3(\hat{x}_t)(\frac{f(\hat{x}_t)}{g(\hat{x}_t)})''' + 6\mu_2 g^2(\hat{x}_t)g'(\hat{x}_t)(\frac{f(\hat{x}_t)}{g(\hat{x}_t)})''$$
$$+ 3\mu_2 g^2(\hat{x}_t)g''(\hat{x}_t)(\frac{f(\hat{x}_t)}{g(\hat{x}_t)})' + 6\mu_2 g(\hat{x}_t)(g'(\hat{x}_t))^2(\frac{f(\hat{x}_t)}{g(\hat{x}_t)})'$$
$$+ \mu_1 g(\hat{x}_t)g''(\hat{x}_t) + \mu_1 (g'(\hat{x}_t))^2)$$
$$- P_t^2 \varphi_\eta^{-1}(h'(\hat{x}_t))^2)dt + P_t^2 h''(\hat{x}_t)\varphi_\eta^{-1}(dz_t - (h(\hat{x}_t) + \frac{1}{2}P_t h''(\hat{x}_t))dt), \qquad (22)$$

where the prime sign in the above expressions indicates derivative with respect to the conditional expectation of the phase variable $x_t$. We state two special cases of equations (21)-(22): (i) For the term $\mu_2 = 0$, the above approximate filtering equations (21)-(22) reduce to the *classical approximate filtering equations* resulting from the second-order approximation (Pugachev and Synstin 1987, p. 452) (ii) the terms $\mu_1 = 2D$, and $\mu_1 = D\tau_{cor}$ of equation (10) with equations (21)-(22) lead to the second-order filtering equations of the OU process-driven stochastic differential system. Hence, equations (21)-(22) describing approximate filtering equations are more general in contrast to the *classical approximate filtering* equations and the *OU process filtering*. The numerical coefficients of equations (21)-(22) are attributed to the successive differentiation with respect to the conditional expectation of the phase variable $x_t$.

## 3. A scalar non-Markovian Duffing system

The theory of this paper is applied to the scalar Duffing system driven by the OU process. Consider a second-order fluctuation equation, i.e.

$$\ddot{x}_t = \alpha x_t + \beta \dot{x}_t - a x_t^3 + g(x_t)\xi_t,$$

where $\alpha < 0$, $\beta \gg 0$, and $\xi_t$ is a Gauss-Markov process. After accomplishing the phase space formulation, we have

$$x_t = x_1, \quad \dot{x}_1 = x_2, \quad \dot{x}_2 = \alpha x_1 + \beta x_2 - a x_1^3 + g(x_1)\xi_t.$$

Consider the contribution to the evolution of the phase variable $x_t = x_1$ coming from the damping term is considerably greater than the inertial term (Hänggi 1995, p. 241), then the term $\alpha x_1 + \beta x_2 - a x_1^3 + g(x_1)\xi_t$ vanishes, i.e.

$$\dot{x}_t = -\frac{\alpha}{\beta}x_t + \frac{a}{\beta}x_t^3 - \frac{g(x_t)}{\beta}\xi_t, \tag{23}$$

The above equation (23) describes a scalar stochastic Duffing system. Equation (23) can be restated as

$$\dot{x}_t = f(x_t) - \frac{g(x_t)}{\beta}\xi_t, \tag{24}$$

where $f(x_t) = -\frac{\alpha}{\beta}x_t + \frac{a}{\beta}x_t^3$. For the above stochastic differential system of equation (24), the approximate filtering equations (21)-(22) can be recast as

$$d\hat{x}_t = (f(\hat{x}_t) - \frac{\mu_2}{2\beta^2}g^2(\hat{x}_t)g'(\hat{x}_t)(\frac{f(\hat{x}_t)}{g(\hat{x}_t)})' - \frac{\mu_2}{2\beta^2}g^3(\hat{x}_t)(\frac{f(\hat{x}_t)}{g(\hat{x}_t)})''$$

$$+ \frac{1}{2}P_t(f''(\hat{x}_t) - \frac{\mu_2}{2\beta^2}(g^3(\hat{x}_t)\frac{d^4}{d\hat{x}_t^4}(\frac{f(\hat{x}_t)}{g(\hat{x}_t)}) + 7g^2(\hat{x}_t)g'(\hat{x}_t)(\frac{f(\hat{x}_t)}{g(\hat{x}_t)})'''$$

$$+ 5g^2(\hat{x}_t)g''(\hat{x}_t)(\frac{f(\hat{x}_t)}{g(\hat{x}_t)})'' + 10g(\hat{x}_t)(g'(\hat{x}_t))^2(\frac{f(\hat{x}_t)}{g(\hat{x}_t)})''$$

$$+ g^2(\hat{x}_t)g'''(\hat{x}_t)(\frac{f(\hat{x}_t)}{g(\hat{x}_t)})' + 2(g'(\hat{x}_t))^3(\frac{f(\hat{x}_t)}{g(\hat{x}_t)})'$$

$$+ 6g(\hat{x}_t)g'(\hat{x}_t)g''(\hat{x}_t)(\frac{f(\hat{x}_t)}{g(\hat{x}_t)})')))dt$$

$$+ P_t h'(\hat{x}_t)\varphi_\eta^{-1}(dz_t - (h(\hat{x}_t) + \frac{1}{2}P_t h''(\hat{x}_t))dt), \tag{25}$$

$$dP_t = (2P_t(f'(\hat{x}_t) - \frac{1}{2}\frac{\mu_2}{\beta^2}(4g^2(\hat{x}_t)g'(\hat{x}_t)(\frac{f(\hat{x}_t)}{g(\hat{x}_t)})'' + g^2(\hat{x}_t)g''(\hat{x}_t)(\frac{f(\hat{x}_t)}{g(\hat{x}_t)})'$$

$$+ 2g(\hat{x}_t)(g'(\hat{x}_t))^2(\frac{f(\hat{x}_t)}{g(\hat{x}_t)})' + g^3(\hat{x}_t)(\frac{f(\hat{x}_t)}{g(\hat{x}_t)})'''))$$

$$+ \frac{\mu_1}{\beta^2}g^2(\hat{x}_t)(1 + 2\frac{\mu_2}{\mu_1}g(\hat{x}_t)(\frac{f(\hat{x}_t)}{g(\hat{x}_t)})')$$

$$+ P_t(\frac{\mu_2}{\beta^2}g^3(\hat{x}_t)(\frac{f(\hat{x}_t)}{g(\hat{x}_t)})''' + 6\frac{\mu_2}{\beta^2}g^2(\hat{x}_t)g'(\hat{x}_t)(\frac{f(\hat{x}_t)}{g(\hat{x}_t)})''$$

$$+ 3\frac{\mu_2}{\beta^2}g^2(\hat{x}_t)g''(\hat{x}_t)(\frac{f(\hat{x}_t)}{g(\hat{x}_t)})' + 6\frac{\mu_2}{\beta^2}g(\hat{x}_t)(g'(\hat{x}_t))^2(\frac{f(\hat{x}_t)}{g(\hat{x}_t)})'$$

$$+ \frac{\mu_1}{\beta^2}g(\hat{x}_t)g''(\hat{x}_t) + \frac{\mu_1}{\beta^2}(g'(\hat{x}_t))^2)$$

$$- P_t^2\varphi_\eta^{-1}(h'(\hat{x}_t))^2)dt + P_t^2 h''(\hat{x}_t)\varphi_\eta^{-1}(dz_t - (h(\hat{x}_t) + \frac{1}{2}P_t h''(\hat{x}_t))dt). \tag{26}$$

Consider the system non-linearity $f(x_t) = -\frac{\alpha}{\beta}x_t + \frac{a}{\beta}x_t^3$, the process noise coefficient $g(x_t) = x_t$, and the input noise process $\xi_t$ is the OU process, i.e. $\mu_1 = 2D$, $\mu_2 = D\tau_{cor}$, see equation (10). As a result of these, we have the following second-order filtering equations for the scalar damped Duffing SDE driven by the OU process:

$$d\hat{x}_t = (-\frac{\alpha}{\beta}\hat{x}_t + \frac{a}{\beta}\hat{x}_t^3 - \frac{2a\mu_2}{\beta^3}\hat{x}_t^3 + \frac{1}{2}P_t(\frac{6a\hat{x}_t}{\beta} - \frac{12a\mu_2\hat{x}_t}{\beta^3}))dt + P_t\varphi_\eta^{-1}(dz_t - \hat{x}_t)dt), \tag{27}$$

$$dP_t = (2P_t(-\frac{\alpha}{\beta} + \frac{3a}{\beta}\hat{x}_t^2 - \frac{6a\mu_2}{\beta^3}\hat{x}_t^2) + \frac{\mu_1\hat{x}_t^2}{\beta^2}(1 + \frac{4\mu_2 a\hat{x}_t^2}{\beta\mu_1}) + \frac{1}{2}P_t(\frac{2\mu_1}{\beta^2} + \frac{48a\mu_2\hat{x}_t^2}{\beta^3}) - P_t^2\varphi_\eta^{-1})dt. \tag{28}$$

Note that the evolution of conditional mean, equation (27), assumes the structure of stochastic differential equation resulting from the stochastic term $dz_t$. On the other hand, the evolution of conditional variance, equation (28), assumes the structure of differential equation because of the linear observation equation. The second-order partial $h''(x_t)$ of the term $h(x_t) = a_t x_t$ of the linear observation equation vanishes. As a result of this, the last term of the right-hand side of equation (26) vanishes. Thus, the evolution of conditional variance becomes a differential equation rather than the SDE. If the state vector is not observable, the structure of the observation equation is non-linear. As a result of this, the evolution of conditional variance assumes the structure of the SDE. To prevent the 'negative' variance in the variance trajectory, the random forcing term is ignored that leads to the concept of modified approximate filters (Jazwinski 1970, p. 344).

## 4. Numerical Simulations

Approximate evolution equations, equations (25)-(26), are hard to evaluate theoretically since the global properties are replaced with the local (Jazwinski 1970, p. 360), i.e. replacing the system non-linearity $a(x_t)$, the process noise coefficient $b(x_t)$ and the measurement non-linearity $h(x_t)$ with their derivatives. Here, we accomplish the numerical testing of the conditional mean and variance evolutions for the scalar stochastic damped Duffing system. The filter must be tested under a variety of conditions (Kushner 1967). For this reason, this paper encompasses four different sets of initial conditions and system parameters for the numerical simulation. For the stochastic problem concerned here, the contribution to the evolution of the state coming from the damping term is considerably greater in contrast to the inertial term, which implies the damping term $\beta \gg 1$. The evolution of the state vector of dynamical systems depends on initial conditions. Random initial conditions lead to the non-zero conditional variance matrix $P_{xx}(0)$ of the state vector at the initial time instant. The non-zero diagonal entries of the conditional variance matrix and the zero off-diagonal entries suggest uncertainties in initial conditions and uncertainties are uncorrelated respectively. As a result of these, the contribution to the evolution of conditional variance for the non-linear 'non-autonomous system' comes from initial variances, system non-linearities. This explains the non-linear non-autonomous system with random initial conditions can be regarded as the stochastic differential system (Karatzas and Shreve 1991, p. 284).

The numerical simulation of the conditional mean and variance trajectories of this paper accounts for random initial conditions since random initial conditions explain the actual physical situation. Furthermore, the filtering equations of this paper involve the term $\varphi_n$ describing the intensity of noise introduced into the observation. The condition $\varphi_n^{-1} = 0$ suggests observations are completely masked by the noise that can be regarded as valueless observations. As a result of this, the filtering density reduces to the prediction density. On the other hand, the term $\varphi_n = 0$ implies noiseless observation and filtering equations will not exist. Interestingly, in mathematical control theory, the deterministic counterpart of the stochastic filtering problem is popularly known as the observer problem involving noiseless observations (Sontag 1998, p. 310). For this reason, the effectiveness of the stochastic filtering equations of this paper, described by equations (27)-(28), is examined under less intensity of noise introduced into the observation as well as larger intensity of noise introduced, i.e. $0 < \varphi_n < \infty$. The first set of initial conditions and system parameters is the following: $D = 5$, $\tau_{cor} = 0.005$, $a = 0.001$, $\beta = 10^4$, $\alpha = -0.001$, $\varphi_n = 10^4$, $P_{xx}(0) = 0.01$. The initial conditions and system parameters considered in this paper are associated with equation (23). This set of initial conditions corresponds to the large intensity of noise introduced into observations as well as relatively less uncertainties in initial conditions.

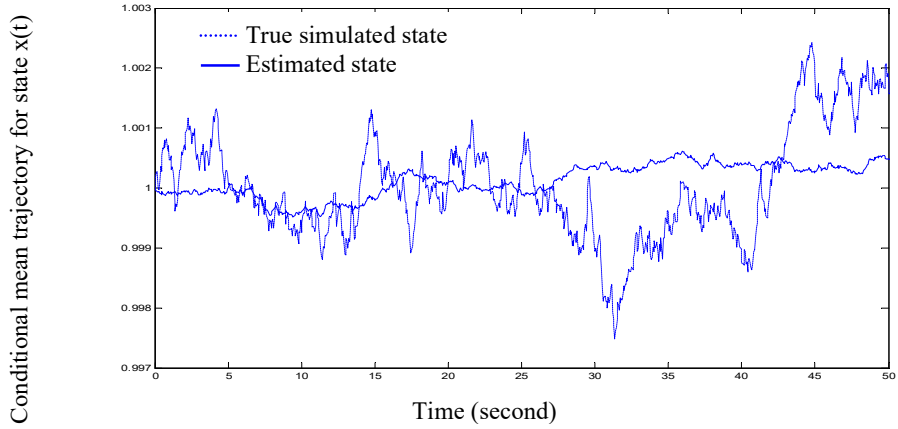

**Fig. 1** a comparison between true and estimated states

Fig. 1 corresponds to the first set of initial conditions and demonstrates that the second-order filter produces bounded state estimates. Here, the second-order filter offers simplified analysis as well as preserves non-linearities via accounting the second-order partial of the system non-linearity, $\frac{1}{2}P_{xx}\frac{\partial^2 a(\hat{x}_t^t)}{\partial \hat{x}_t^2}$. Here we explain 'the filtering equations of this paper preserve perturbation effects felt by the dynamical system'. The stochastic evolution of conditional variance retains the diffusion coefficient $b^2(\hat{x}_t)$ as well as its second-order partial $\frac{1}{2}P_{xx}\frac{\partial^2 b^2(\hat{x}_t)}{\partial \hat{x}_t^2}$. Moreover, the variance term is accounted for in the stochastic evolution of conditional mean. Because of the presence of the variance term in the conditional mean evolution, the second-order filtering accounts for random initial conditions as well as random perturbations.

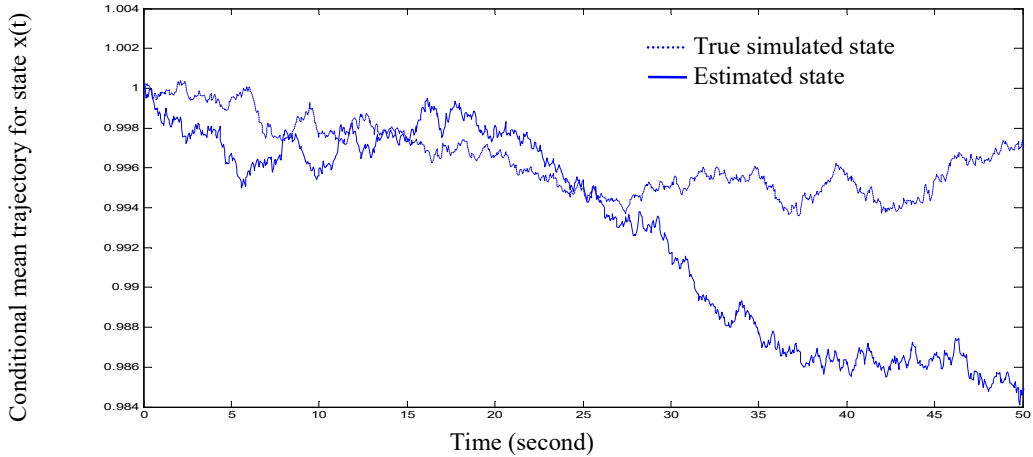

**Fig. 2** a comparison between true and estimated states

Fig. 2 corresponds to the *second* set of initial conditions, i.e. $a = 0.001$, $\beta = 10^4$, $\alpha = -0.001$, $D = 5$, $\tau_{cor} = 0.001$, $\varphi_n = 10^4$, $P_{xx}(0) = 0.1$. Here, the difference between estimated state trajectory and simulated state becomes quite larger, which is attributed to large measurement noise (Germani *et al.* 2007; Sharma 2009) as well as relatively larger value of the variance of the phase variable at the initial time instant.

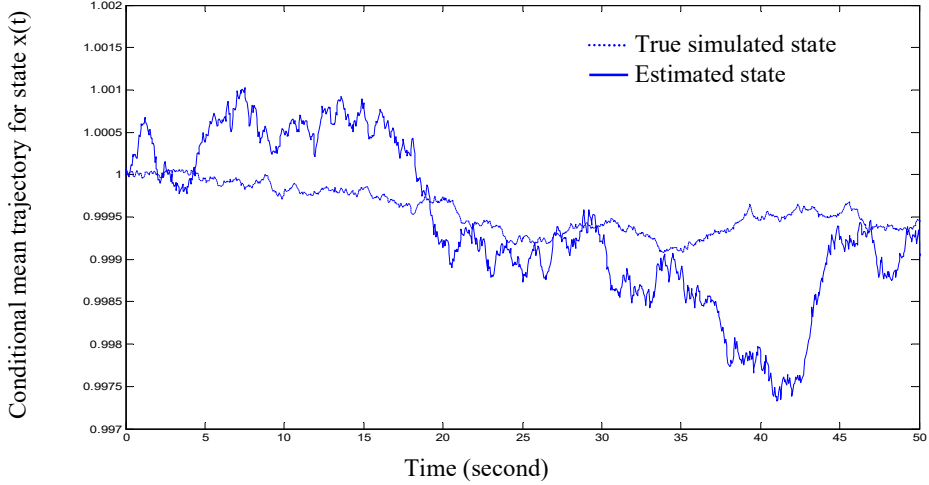

**Fig. 3** a comparison between true and estimated states

Fig. 3 corresponds to the *third* set of initial conditions, i.e. $a = 0.001$, $\beta = 10^3$, $\alpha = -0.001$, $D = 5$, $\tau_{cor} = 0.005$, $\varphi_n = 10^3$, $P_{xx}(0) = 0.1$. The third set of initial conditions is different from the first two sets of initial conditions in the sense that it involves relatively less intensity $\varphi_n$ of the observation noise as well as the smaller value of the damping term $\beta$. The estimated trajectory of Fig. 3 results by simulating the second-order filtering equations, equations (27)-(28) of this paper. The numerical simulation demonstrated in Fig. 3 suggests the efficacy of the second-order filtering equations of this paper for the non-linear state estimation problem since the state estimate does not diverge (Jazwinski 1970, p. 360).

Furthermore, Fig. 4 corresponds to the *fourth* set of initial conditions, i.e. $a = 0.001$, $\beta = 10^3$, $\alpha = -0.001$, $D = 5$, $\tau_{cor} = 0.005$, $\varphi_n = 10^3$, $P_{xx}(0) = 0.01$. This set of initial conditions involves relatively less uncertainties in initial conditions in contrast to the third set. A similar conclusion can be inferred from Fig. 4, which results from the second-order filtering equations of the paper by considering the fourth set of initial conditions.

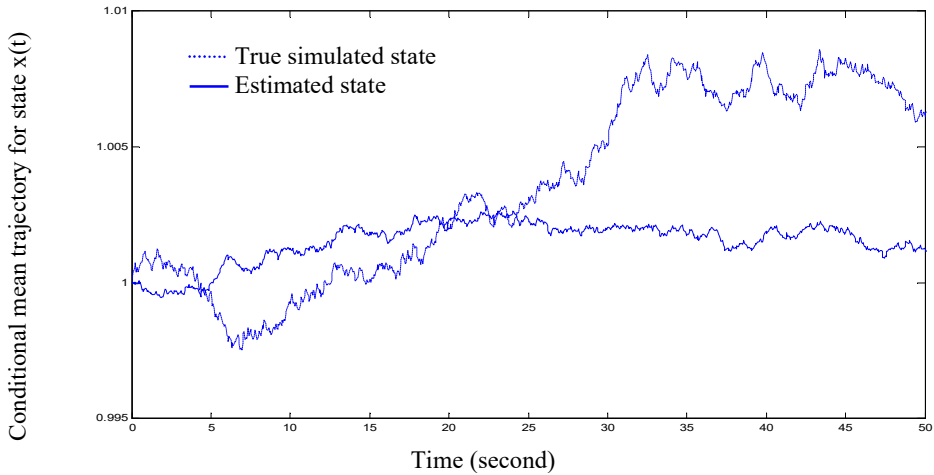

**Fig. 4** a comparison between true and estimated states

### 5. Conclusion

This paper contributes towards the development of a filtering theory for a weakly coloured noise driven non-Markovian stochastic differential system, where the observations are masked by the Brownian noise. The theory of this paper is developed using the standard theory of stochastic processes, i.e. the notion of stochastic equivalence coupled with the filtering density evolution equation. The filtering theory of this paper can be regarded as an extension of 'classical filtering theory'. This paper also contributes to deriving a simpler proof of the Fokker-Planck equation for the OU process-driven stochastic differential system in contrast to the functional calculus approach. The proposed filtering theory is successfully applied to a scalar Duffing system driven by weakly coloured noise. Numerical simulation is carried with different varieties of condition, i.e. changing the intensity, damping term and uncertainties. The results obtained from the proposed filtering theory confirm the effectiveness and usefulness of the method of the paper to estimate the states of a non-Markovian stochastic differential system driven by the OU process. The results of this paper will be useful for filtering and analyzing stochastic problems in which dynamical systems are influenced by weakly coloured noise environment and the tracking station is masked by the Brownian noise. The results of this paper can be extended for the vector case by introducing the component-wise analysis.

**Appendix A:** *The Fokker-Planck equation*

The Fokker-Planck equation is the evolution of conditional probability density for given initial states for a Markov process. The evolution of the conditional probability density of the non-Markov process is known as the stochastic equation. The stochastic equation for a non-Markov process, which satisfies $\dot{x}_t = f(x_t) + g(x_t)\xi_t$, where the input process $\xi_t$ has relatively smaller correlation time, reduces to the Fokker-Planck equation. Consider the input process $\xi_t$ is OU process with smaller correlation time, the Fokker-Planck equation can be derived using the functional calculus (Hänggi 1995). Here, we state an alternative and a brief approach to arrive at the Fokker-Planck equation for the OU process-driven stochastic differential system. The standard structure of the Fokker-Planck equation is given in Stratonovich (1963, p. 62), i.e. $\dot{p}(x) = -\frac{\partial}{\partial x}k_1(x)p(x) + \frac{1}{2}\frac{\partial^2 k_2(x)p(x)}{\partial x^2}$, and consider

$$k_1(x) = M(x) + \frac{1}{4}\frac{dk(x)}{dx}, \tag{A1}$$

$$k_2(x) = k(x), \tag{A2}$$

a simple calculation will show that the Fokker-Planck equation can be recast as

$$\dot{p}(x) = -\frac{\partial}{\partial x}M(x)p(x) + \frac{1}{4}\frac{\partial}{\partial x}(k(x)\frac{\partial p}{\partial x} + \frac{\partial k(x)p(x)}{\partial x}). \tag{A3}$$

The terms $k_1(x)$ and $k_2(x)$ for the OU process-driven stochastic differential system are

$$k_1(x) = f + Dgg' + D\tau_{cor}g^2 g'(\frac{f}{g})', \tag{A4}$$

$$k_2(x) = 2Dg^2 + 2D\tau_{cor}g^3(\frac{f}{g})'. \tag{A5}$$

Equations (4)-(5) in combination with equation (10) result equations (A4)-(A5). Equations (A4)-(A5) in conjunction with equations (A1)-(A2) lead to

$$k(x) = 2Dg^2 + 2D\tau_{cor}g^3(\frac{f}{g})',$$

$$\frac{1}{4}\frac{\partial k(x)}{\partial x} = Dgg' + \frac{3}{2}D\tau_{cor}g^2 g'(\frac{f}{g})' + \frac{1}{2}D\tau_{cor}g^3(\frac{f}{g})'',$$

$$M(x) = f - \frac{1}{2}D\tau_{cor}g^2 g'(\frac{f}{g})' - \frac{1}{2}D\tau_{cor}g^3(\frac{f}{g})''. \tag{A6}$$

Equation (A3) in combination with equation (A6) leads to

$$\dot{p}(x) = -\frac{\partial}{\partial x}(f - \frac{D\tau_{cor}}{2}g^2 g'(\frac{f}{g})' - \frac{1}{2}D\tau_{cor}g^3(\frac{f}{g})'')p + \frac{1}{4}\frac{\partial}{\partial x}((2Dg^2 + 2D\tau_{cor}g^3(\frac{f}{g})')\frac{\partial p}{\partial x}$$

$$+ \frac{\partial}{\partial x}(2Dg^2 + 2D\tau_{cor}g^3(\frac{f}{g})')p).$$

Further simplification leads to

$$\dot{p}(x) = -\frac{\partial}{\partial x}fp + \frac{D\tau_{cor}}{2}\frac{\partial}{\partial x}g^2 g'(\frac{f}{g})'p + \frac{1}{2}D\tau_{cor}\frac{\partial}{\partial x}g^3(\frac{f}{g})''p$$

$$+\frac{1}{4}\frac{\partial}{\partial x}((2Dg^2+2D\tau_{cor}g^3(\frac{f}{g})')\frac{\partial p}{\partial x})+\frac{1}{4}\frac{\partial}{\partial x}((2Dg^2+2D\tau_{cor}g^3(\frac{f}{g})')\frac{\partial p}{\partial x})$$

$$+\frac{1}{4}\frac{\partial}{\partial x}(4Dgg'+2D\tau_{cor}g^3(\frac{f}{g})''+6D\tau_{cor}g^2g'(\frac{f}{g})')p.$$

A re-arrangement of the right-hand side terms of the above equation results

$$\dot{p}(x)=-\frac{\partial}{\partial x}fp+\frac{\partial}{\partial x}((Dg^2+D\tau_{cor}g^3(\frac{f}{g})')\frac{\partial p}{\partial x})+\frac{\partial}{\partial x}(Dgg'+D\tau_{cor}g^3(\frac{f}{g})''+2D\tau_{cor}g^2g'(\frac{f}{g})')p,$$

Furthermore,

$$\dot{p}(x)=-\frac{\partial}{\partial x}fp+D\frac{\partial}{\partial x}g((g+\tau_{cor}g^2(\frac{f}{g})')\frac{\partial p}{\partial x})+D\frac{\partial}{\partial x}g(g'+\tau_{cor}g^2(\frac{f}{g})''+2\tau_{cor}gg'(\frac{f}{g})')p,$$

after combining the last two terms of the above equation, we have

$$\dot{p}(x)=-\frac{\partial}{\partial x}fp+D\frac{\partial}{\partial x}(g\frac{\partial}{\partial x}(g+\tau_{cor}g^2(\frac{f}{g})')p)=-\frac{\partial}{\partial x}fp+D\frac{\partial}{\partial x}(g\frac{\partial}{\partial x}g(1+\tau_{cor}g(\frac{f}{g})')p), \quad (A7)$$

equation (A7) can be regarded as the Fokker-Planck equation for the OU process-driven stochastic differential system, where the OU process is a weakly coloured noise process.